\documentclass[12pt]{amsart}
\usepackage{amssymb}
\usepackage{times}
\usepackage{graphicx}

\theoremstyle{plain}

\numberwithin{equation}{section}

\newcommand{\calD}{\mathcal{D}}

\newcommand{\calH}{\mathcal{H}}

\newcommand{\calS}{\mathcal{S}}

\begin{document}
\title [The moduli space of Hessian quartic surfaces] {The moduli space of Hessian quartic surfaces and automorphic forms}
\author{Shigeyuki Kond{$\bar{\rm o}$}}
\address{Graduate School of Mathematics, Nagoya University, Nagoya,
464-8602, Japan}
\email{kondo@math.nagoya-u.ac.jp}
\thanks{Research of the author is partially supported by
Grant-in-Aid for Scientific Research (S), No 22224001, (S), No 19104001}

\begin{abstract}
We shall show the existence of 15 automorphic forms of weight 8 on the moduli space of marked Hessian quartic surfaces of cubic surfaces. These 15 automorphic forms correspond to $(\lambda_i - \lambda_j)(\lambda_k-\lambda_l)$ where $\lambda_i, \lambda_j, \lambda_k, \lambda_l$ 
$(1\leq i,j,k,l \leq 5)$ are the coefficients of the Sylvester form of a general cubic surface and all $i,j,k,l$ are distinct.
\end{abstract}
\maketitle

\section{Introduction}

The purpose of this note is to give an example of automorphic forms on
the moduli space of Hessian quartic surfaces which can be interpreted in terms of invariants of cubic surfaces.
Let $S$ be a smooth cubic surface defined by a homogeneous polynomial $F(z_0,z_1,z_2,z_3)$ of degree 3.  
Then the hessian of $F$, if it is not identically zero, defines a quartic surface $H$ called the {\it Hessian quartic surface} of $S$.  
To study Hessian quartic surfaces, it is convenient to use the Sylvester form of $S$.
It is classically known (cf. Segre \cite{Se}, Chap. IV) that a general cubic surface can be written by the Sylvester form: 
$$\lambda_1x_1^3 + \cdots + \lambda_5x_5^3 = 0, \quad x_1+\cdots + x_5 =0.$$
A general cubic surface defined by the Sylvester
form is uniquely determined by 
$$\lambda = (\lambda_1: \cdots : \lambda_5)\in {\bf P}^4$$ 
up to permutation of $\lambda_i$.  
By using the theory of periods of $K3$ surfaces (Piatetskii-Shapiro, Shafarevich \cite{PS}),  
one can describe the moduli space of Hessian quartic surfaces as 
an arithmetic quotient of the 4-dimensional bounded symmetric domain of type IV. 
Koike \cite{Koi} gave an ${\mathfrak S}_5$-equivariant birational map from ${\bf P}^4$
to the moduli space of marked Hessian quartic surfaces.  On the other hand, 
Dardanelli and van Geemen \cite{DvG} studied the transcendental lattices of the Hessian quartic surfaces of cubic surfaces with a node, 
with an Eckardt point or without a Sylvester form.
In the moduli space of marked Hessian quartic surfaces,
the Heegner divisor corresponding to cubic surfaces with an Eckardt point consists of 10 irreducible components.  
A cubic surface defined by a Sylvester form has an Eckardt point iff $\lambda_i = \lambda_j$ for some $i\not=j$.
Thus 10 components defined by $\lambda_i = \lambda_j$ in ${\bf P}^4$ bijectively correspond to 10 components of the
above Heegner divisor in the moduli space of marked Hessian quartic surfaces (Lemma \ref{norm1}).  
The minimal model of a general Hessian quartic surface has a canonical fixed point free involution 
(for example, see Dolgachev, Keum \cite{DK}), and hence it is the covering $K3$ surface of an Enriques surface.
Thus the moduli space of cubic surfaces is birational to the moduli space of Enriques surfaces whose covering $K3$ surfaces are Hessian
quartic surfaces.
In the paper \cite{Kon}, the author gave $3^3\cdot 5\cdot 17 \cdot 31$ holomorphic automorphic forms $F_V$ of weight 4 with known zeros
on the moduli space of marked Enriques surfaces 
by using Borcherds theory of automorphic forms \cite{B}.  

In this note, 
by first dividing $F_V$ by a product of suitable linear forms and restricting it to
the locus of marked Enriques surfaces corresponding to
marked Hessian quartic surfaces, we have 15 automorphic forms of weight 8. 
We shall show that under the birational map 15 automorphic forms correspond to 
$(\lambda_i - \lambda_j)(\lambda_k-\lambda_l)$ where all $i, j, k, l$ are distinct  (Theorem \ref{rest2}).

The Hessian quartic surface is an example of the Cayley symmetroid given by the zero of a $4 \times 4$ symmetrical
determinant whose entries are the linear form (Dardanelli, van Geemen \cite{DvG}, \S 1.6).  
A general Cayley symmetroid is the covering $K3$ surface of a nodal Enriques surface, namely an Enriques surface containing a smooth rational curve (Cossec \cite{Cos}).
The period domain of nodal Enriques surfaces is the $9$-dimensional bounded symmetric domain ${\calD}$ of type IV which is naturally embedded in ${\calD}(L_-)$. 
By the similar way as in the case of Hessian quartic surfaces, we have $3^3\cdot 5\cdot 17$ holomorphic automorphic forms on ${\calD}$. It would be interesting to study a relation between these automorphic forms and the geometry of Cayley symmetroids given in Coble \cite{C}, Chapters V, VI.

Finally we mention the related works.  
In the paper \cite{ACT}, Allcock, Carlson, Toledo showed that the moduli space of cubic surfaces can be described as an arithmetic 
quotient of the 4-dimensional complex ball by considering the period of the  intermediate Jacobian of the triple cover of ${\bf P}^3$ branched along a cubic surface.  Later Dolgachev, van Geemen and the author \cite{DvGK} 
gave the same description of the moduli space of cubic surfaces by using the theory of periods of $K3$ surfaces.
By using this description of the moduli space and Borcherds' theory of automorphic forms \cite{B}, 
Allcock, Freitag \cite{AF} studied the moduli space of marked cubic surfaces.  They constructed automorphic forms corresponding to
Cayley's cross ratios for cubic surfaces.

\smallskip

The author would like to thank Igor Dolgachev for useful conversations.

\section{Preliminaries}\label{}

A {\it lattice} $(L, \langle, \rangle)$ is a pair of a free ${\bf Z}$-module $L$ of rank $r$ and
a non-degenerate symmetric integral bilinear form $\langle, \rangle : L \times L \to {\bf Z}$. 
For simplicity we omit $\langle, \rangle$ if there are no confusions.  For $x \in L\otimes {\bf Q}$, we call $x^2 =\langle x,x\rangle$ the {\it norm} of $x$.
For a lattice $(L,\langle,\rangle)$ and an integer
$m$, we denote by $L(m)$ the lattice $(L, m\langle, \rangle)$.  We denote by $U$ the even unimodular lattice of signature $(1,1)$, 
and by $A_m, \ D_n$ or $\ E_k$ the even {\it negative} definite lattice defined by
the Cartan matrix of type $A_m, \ D_n$ or $\ E_k$ respectively.  For an integer $m$, we denote by $\langle m\rangle$ the lattice of rank 1
generated by a vector with norm $m$.  We denote by $L\oplus M$ the orthogonal direct sum of lattices $L$ and $M$.

Let $L$ be an even lattice and let $L^* ={\rm Hom}(L,{\bf Z})$.  We denote by $A_L$ the quotient
$L^*/L$ and define a map
$$q_L : A_L \to {\bf Q}/2{\bf Z}$$
by $q_L(x+L) = \langle x, x\rangle\ {\rm mod}\ 2{\bf Z}$.  We call $q_L$ the {\it discriminant quadratic form} of $L$.  
We denote by $u$ or $v$ the discriminant quadratic form of $U(2)$ or $D_4$ respectively.

Let ${\rm O}(L)$ be the orthogonal group of $L$, that is, the group of isomorphisms of $L$ preserving the bilinear form.
Similarly ${\rm O}(q_L)$ denotes the group of isomorphisms of $A_L$ preserving $q_L$.
There is a natural map
$${\rm O}(L) \to {\rm O}(q_L)$$
whose kernel is denoted by $\tilde{\rm O}(L)$.

\section{The Hessians of cubic surfaces and Enriques surfaces}\label{sec:2}

Let $S$ be a smooth cubic surface defined by a homogeneous polynomial $F(z_0,z_1,z_2,z_3)$ of degree 3.  
Then the hessian of $F$, if it is not identically zero, defines a quartic surface $H$ called the {\t Hessian quartic surface} of $S$.  
It is classically known that a general cubic surfaces $S$ can be written by the {\it Sylvester form}
\begin{equation}\label{sylvester}
\lambda_1x_1^3 + \cdots + \lambda_5x_5^3 = 0, \quad x_1+\cdots + x_5 =0
\end{equation}
where $x_1,..., x_5$ are linear forms in $z_0,z_1,z_2,z_3$ each four of them are linearly independent and $\lambda_i \in {\bf C}^*$.
The forms $x_1,..., x_5$ are uniquely determined by $F$ up to permutation and multiplication by a common non-zero scalar,
and $\lambda_1,..., \lambda_5$ are uniquely determine by $F$ and $x_i$.  Thus a general cubic surface defined by the Sylvester
form is now determined by 
$$\lambda = (\lambda_1: \cdots : \lambda_5)\in {\bf P}^4$$ 
up to permutations
of $\lambda_i$ (Segre \cite{Se}, Chap. IV).

For a cubic surface defined by the Sylvester form, the corresponding Hessian quartic surface $H$ is given by
\begin{equation}\label{hessian}
{1 \over \lambda_1x_1} + \cdots + {1\over \lambda_5x_5} = 0, \quad x_1+\cdots + x_5 =0.
\end{equation}
The Hessian quartic surface $H$ has 10 nodes $p_{ijk}$ defined by $x_i = x_j = x_k =0$, and contains 10 lines $l_{mn}$ defined by $x_m = x_n =0$.
It is known  (Segre \cite{Se}, Chap. IV) that $H$ has no other singular points if and only if 
\begin{equation}\label{nodal}
\Delta_{\rm sing}(\lambda) = \sum_{i=1}^5 {1\over \pm \sqrt{\lambda_i}} \not= 0.
\end{equation}
We denote by $X$ the minimal resolution of $H$
which is a $K3$ surface with 20 smooth rational curves, that is, exceptional curves $E_{ijk}$ over 10 nodes $p_{ijk}$ 
and strict transforms $L_{mn}$ of 10 lines $l_{mn}$.  The curve $E_{ijk}$ meets exactly three curves $L_{ij}$, $L_{ik}$ and $L_{jk}$, and
conversely $L_{ij}$ meets exactly three curves $E_{ijk} \ (k\not=i,j)$.
Thus we have two sets $\{E_{ijk}\}$, $\{L_{mn}\}$ of smooth rational curves on $X$ each of which consists of 10 disjoint curves, and
each curve in one set meets exactly three curves in the other set.

The birational involution defined by
\begin{equation}\label{involution}
 (x_1:\cdots : x_5) \to ({1 \over \lambda_1x_1} : \cdots : {1\over \lambda_5x_5})
\end{equation}
induces a fixed point free involution $\sigma$ of $X$, and hence the quotient $Y = X/\langle \sigma \rangle$ is an Enriques surface
(Dolgachev, Keum \cite{DK}).
The involution $\sigma$ switches nodal curves $E_{ijk}$ and $L_{mn}$ where $\{i,j,k,m,n\} = \{ 1,2,3,4,5\}$.
We denote by $\bar{L}_{ij}$ the image of $L_{ij}$ or $E_{kmn}$ on $Y$.  
The curve $\bar{L}_{ij}$ meets exactly three curves $\bar{L}_{km}$, $\bar{L}_{kn}$ and
$\bar{L}_{mn}$.   We can easily see that the dual graph of ten nodal curves $\{ \bar{L}_{ij}\}$ is 
isomorphic to the Petersen graph whose automorphism group is the symmetry group ${\mathfrak S}_5$ of degree $5$.

Denote by
$\pi : X \to Y$ the natural projection and
let $L$ be the lattice $H^2(X, {\bf Z})$ which is the even unimodular lattice of signature $(3,19)$.
Define
\begin{equation}\label{invar}
L_{\pm} = \{ x \in L \ | \ \sigma^*(x) = \pm x \}.
\end{equation}
It is known that $L_+ \cong \pi^*(H^2(Y,{\bf Z})) \cong U(2) \oplus E_8(2)$ and $L_-\cong U\oplus U(2)\oplus E_8(2)$ (Barth, Peters \cite{BP}).

The 20 curves $\{ E_{ijk}, L_{mn}\}$ generate a sublattice $N$ of signature $(1,15)$ in the Picard lattice $S_X$ of $X$.  
Let $M$ be the orthogonal complement of $N$ in $L$.  It is known that $M$ is isomorphic to $U\oplus U(2) \oplus A_2(2)$
(Dolgachev, Keum \cite{DK}).
Let $R$ be the orthogonal complement of $M$ in $L_-$.  Obviously $R$ is a negative definite lattice of rank $6$.

\subsection{Lemma}\label{root-inv}
{\it $R$ is isomorphic to $E_6(2)$.}

\begin{proof}
For any smooth rational curve $C$ on $Y$, $\pi^*(C)$ is the disjoint union of two smooth rational curves.
The difference of these two curves is a vector of norm $(-4)$ contained in $R$. For example, 
$E_{123}-L_{45}$, $E_{145}-L_{23}$, $E_{235}-L_{14}$, $E_{345}-L_{12}$, $E_{125}-L_{34}$, and $E_{245}-L_{13}$ generate
a lattice isomorphic to $E_6(2)$ in $R$.  
By comparing $A_{L_-} = ({\bf Z}/2{\bf Z})^{10}$ and 
$A_{E_6(2)\oplus M} = ({\bf Z}/2{\bf Z})^{10}\oplus ({\bf Z}/3{\bf Z})^2$, 
we can conclude that $E_6(2)$ is the orthogonal complement of $M$ in $L_-$.
\end{proof}

\subsection{Remark}\label{}
{\it Nikulin {\rm \cite{N2}} introduced the notion of {\it root invariant} $(R', K)$ of each Enriques surface consisting of a root lattice $R'$ and 
a finite subgroup $K$ of $R'/2R'$.
The $R'(2)$ is generated by the differences $C-\pi^*(C)$ of all smooth rational curves $C$ on the Enriques surface.
In case that the covering $K3$ surface of $Y$ is a Hessian quartic surface $X$, the generic $Y$ has the root invariant $(E_6, \{0\})$.}

\medskip

Let $S_X$ be the Picard lattice of $X$.  The orthogonal complement of $S_X$ in 
$H^2(X,{\bf Z})$, denoted by $T_X$, is called the {transcendental lattice} of $X$.
For a generic cubic surface $S$, $N$ (resp. $M$) coincides with $S_X$ (resp. $T_X$) (Dolgachev, Keum \cite{DK}).

Recall that a smooth cubic surface $S$ has 45 tritangent planes each of which consists of three lines.
If three coplanar lines meet at one point, the intersection point is called an {\it Eckardt point}.
A smooth cubic surface $S$ given by the Sylvester form (\ref{sylvester}) has an Eckardt point if and only if $\lambda_i = \lambda_j$
(Segre \cite{Se}, Chap. IV).
If $S$ has an Eckardt point, for example $\lambda_i = \lambda_j$, the tritangent plane is given by $x_i+x_j =0$ and 
the Eckardt point on $S$ is $p_{kmn}$
where $\{ i,j,k,m,n\} = \{1,2,3,4,5\}$.
In this case, the plane section defined by $x_i+x_j=0$ on the Hessian quartic surface consists of $2l_{ij}$ and two lines through the node
$p_{kmn}$ (see \cite{DvG}, \S 2.1, 2.2).  
The strict transforms of the two lines to $X$ are two
disjoint smooth rational curves $N_{ij}^+, N_{ij}^-$.  The involution $\sigma$ switches $N_{ij}^+$ and  $N_{ij}^-$.
The curves $N_{ij}^{\pm}$ meet $L_{ij}$ and $E_{kmn}$ with multiplicity 1 and disjoint with other 18 curves.
Thus we have the following Lemma.

\subsection{Lemma}\label{eckardt}
{\it If a smooth cubic surface has an Eckardt point corresponding to $\lambda_i=\lambda_j$, 
then $X$ contains new two smooth rational curves $N_{ij}^+$ and $N_{ij}^-$.
The class of $N_{ij}^+ - N_{ij}^-$ is a $(-4)$-vector contained in $M$.}

\begin{proof}
Note that $N_{ij}^+ - N_{ij}^-$ is perpendicular to 20 smooth rational curves
$\{ L_{ij}, E_{kmn}\}$.  Since 20 curves $\{ L_{ij}, E_{kmn}\}$ generate $N$, we have $N_{ij}^+ - N_{ij}^- \in N^{\perp}=M$. 
\end{proof}

\subsection{Remark}\label{clebsch}
{\it If all $\lambda_i =1$, the cubic surface is called the Clebsch diagonal cubic surface which has 
$10$ Eckardt points
{\rm (}see {\rm \cite{DvG}, Lemma 2.2)}.
The corresponding Enriques surface contains $20$ smooth rational curves and the symmetry group ${\mathfrak S}_5$ of degree $5$ acts
on the Enriques surface as automorphisms.  This Enriques surface is one of Enriques surfaces with a finite group of automorphisms
classified in {\rm \cite{Kon1}} {\rm (}see {\rm \cite{Kon1}, Example VI, \cite{DvG}, \S 2.3)}.}

\section{Discriminant quadratic form}\label{discriminant}

First recall that $N$ and $M$ are primitive sublattices of the even unimodular lattice $L = H^2(X,{\bf Z})$ with $M=N^{\perp}$.
It follows from Nikulin \cite{N1}, Corollary 1.6.2 that $q_M \cong -q_N$.  
By an elementary calculation, $A_M \cong ({\bf Z}/2{\bf Z})^4 \oplus {\bf Z}/3{\bf Z}$ and
the restriction $(q_M)_2$ of $q_M$ to the 2-Sylow subgroup of $A_M$ is isomorphic to 
$u\oplus v$.  We can consider $(q_M)_2$ a 4-dimensional quadratic form over ${\bf F}_2$.  It is well known that 
the group of automorphisms of the quadratic form $(q_M)_2$ is
isomorphic to the symmetry group ${\mathfrak S}_5$ of degree 5 (\cite{CC}, page 2), and hence ${\rm O}(q_M)$ is isomorphic to
${\mathfrak S}_5 \times {\bf Z}/2{\bf Z}$ where ${\bf Z}/2{\bf Z}$ is the involution of  ${\bf Z}/3{\bf Z}$.
It is easily see that $(q_M)_2$ contains $2(2^2-1)$ isotropic vectors
and $2(2^2+1)$ non-isotropic vectors.  For $a\in A_M$ we denote by $|a|$ the order of $a$.
We can easily see the following lemma. 

\subsection{Lemma}\label{discri}
{\it The number of vectors $a$ in $A_M$ with 

$(|a|, q_M(a)) = (0,0), (2,0), (2,1), (3, -4/3), (6, -1/3)$ or $(6,-4/3)$

\noindent
is 
$1, 5, 10, 2, 20$ or $10$, respectively.
}

\medskip
\noindent
It follows from Lemma \ref{discri} that $A_N$ ($\cong A_M$) contains 10 vectors with norm 1.
In the following we shall study a geometric meaning of these 10 vectors.
Recall that the Enriques surface $Y$ contains 10 smooth rational curves $\bar{L}_{ij}$.
If we fix $\bar{L}_{ij}$ one of them, then there are 6 smooth rational curves perpendicular to $\bar{L}_{ij}$ which form
a singular fiber of type $I_6$ of an elliptic fibration.  We denote this elliptic fibration by $|\bar{F}_{ij}|$.
For example, if we take $\bar{L}_{12}$, then the class
$$\bar{F}_{12} = \bar{L}_{13} + \bar{L}_{24}+\bar{L}_{15} + \bar{L}_{23} + \bar{L}_{14} + \bar{L}_{25} = 
\bar{E}_{245} + \bar{E}_{135}+\bar{E}_{234} + \bar{E}_{145} + \bar{E}_{235} + \bar{E}_{134}$$
defines an elliptic fibration on $Y$, and $\bar{L}_{12}$ is a component of an another singular fiber of the fibration.
Then $\pi^*(\bar{F}_{ij}) = 2F_{ij}$ and the class
$$F_{ij} = L_{13}+E_{135}+L_{15}+E_{145}+L_{14}+E_{134} = E_{245}+L_{24}+E_{234}+L_{23}+E_{235}+L_{25}$$ 
defines an elliptic fibration on $X$.
We can easily see that 
$$\alpha_{ij} = {1\over 2}(F_{ij} -  L_{ij} - E_{kmn}), \quad \{i,j,k,m,n\} = \{1,2,3,4,5\}$$
has an integral intersection number with any curve from 20 curves $\{ L_{ij}, E_{kmn}\}$.  
Since 20 curves $\{ L_{ij}, E_{kmn}\}$ generate $N$, $\alpha_{ij}$ is contained in $N^*$.
Obviously $q_N(\alpha_{ij})=1$.

In the Sylvester form (\ref{sylvester}), if $\lambda_i =\lambda_j$, then $S$ has an Eckardt point $p_{kmn}$.
We have two new smooth rational curves $N_{ij}^+, \ N_{ij}^-$ meeting $L_{ij}$ and $E_{kmn}$ (Lemma \ref{eckardt}).
Recall that the involution $\sigma$ switches $N_{ij}^+$ and $N_{ij}^-$.  Let $\bar{N}_{ij}$ be the image of $N_{ij}^+$ on $Y$.
Then $\bar{N}_{ij} + \bar{L}_{ij}$ is a singular fiber of $|\bar{F}_{ij}|$ of type $I_2$.  
Note that this is a multiple fiber because $N_{ij}^+ + N_{ij}^- + L_{ij} + E_{kmn}$ is a singular fiber of type ${\rm I}_4$ of
the elliptic fibration $| F_{ij}|$.
It follows that  
$$\alpha_{ij} = {1\over 2}(N_{ij}^+ + N_{ij}^-).$$
The difference 
$$\beta_{ij} = {1\over 2}(N_{ij}^+ - N_{ij}^-)$$
defines a vector in $M^*$ with $q_M(\beta_{ij})=1$ (Lemma \ref{eckardt}).  The condition 
$$\alpha_{ij} + \beta_{ij} \in L=H^2(X,{\bf Z})$$ 
gives a bijective correspondence
between the set of vectors with norm 1 in $A_M$ and that in $A_N$.

Moreover if $\lambda_1 =\lambda_2$, $\lambda_3 =\lambda_4$, but other coefficients are different, then
$S$ has exactly two Eckardt points.
In this case, $\beta_{12}$ is perpendicular to $\beta_{34}$ (see Dardanelli, van Geemen \cite{DvG}, Lemma 2.2, the case $k =2$).
Thus we have the following Lemma.

\subsection{Lemma}\label{norm1}
{\it There is a bijective correspondence between the $10$ conditions $\lambda_i = \lambda_j$ between $\{\lambda_i\}$ and
the set of vectors in $A_M$ with norm $1$.  Moreover if $S$ has exactly two Eckardt points corresponding to $\lambda_i =\lambda_j$ and
$\lambda_k =\lambda_m$ where all $i,j,k,m$ are distinct, then $\beta_{ij}$ is perpendicular to $\beta_{km}$.
}

\section{Periods and Heegner divisors}\label{periods}

In this section, we recall the period domain and Heegner divisors for Enriques surfaces and Hessian quartic surfaces.
First we consider the case of Enriques surfaces.
Define
\begin{equation}\label{period1}
{\calD}(L_-) = \{ [\omega] \in {\bf P}(L_-\otimes {\bf C}) \ | \ \langle \omega, \omega \rangle = 0, \ \langle \omega, \bar{\omega}\rangle > 0\}
\end{equation}
which is a disjoint union of two copies of the 10-dimensional bounded symmetric domain of type IV.

The discriminant quadratic form $(A_{L_-}, q_{L_-})$ is the orthogonal direct sum of five copies of $u$.
We consider the orthogonal group ${\rm O}(L_-)$ of $L_-$ and denote by $\tilde{\rm O}(L_-)$ the kernel of the map
$${\rm O}(L_-) \to {\rm O}(q_{L_-}).$$
Then ${\rm O}(L_-)/\tilde{\rm O}(L_-) \cong {\rm O}(q_{L_-}) \cong {\rm O}^+(10, {\bf F}_2)$ (Barth, Peters \cite{BP}).

For a vector $r \in L_-$ with a negative norm, we put
$$r^{\perp} = \{ [\omega] \in {\calD}(L_-) \ | \ \langle \omega, r \rangle = 0\}.$$
Let $a \in A_{L_-}$ be a non-isotropic vector, that is, $q_{L_-}(a) = 1$.
We define Heegner divisors $\tilde{\calH}$ and $\tilde{\calH}_a$ by
$$\tilde{\calH} = \sum_r \ r^{\perp}, \quad \tilde{\calH}_a = \sum_t \ t^{\perp}$$
where $r$ moves over the set of all $(-2)$-vectors in $L_-$ and $t$ moves over the set of all $(-4)$-vectors in $L_-$ satisfying ${t \over 2} \ {\rm mod} \ L_- = a$.  
It is known that ${\calD}(L_-)\setminus \tilde{\calH}$ is the period domain of Enriques surfaces.
The quotient $({\calD}(L_-)\setminus \tilde{\calH})/{\rm O}(L_-)$ (resp. $({\calD}(L_-)\setminus \tilde{\calH})/\tilde{\rm O}(L_-)$) is the moduli space of Enriques surfaces (resp. 
the moduli space of marked Enriques surfaces).

The case of Hessian quartic surfaces is similar.  First define
\begin{equation}\label{period2}
{\calD}(M) = \{ [\omega] \in {\bf P}(M\otimes {\bf C}) \ | \ \langle \omega, \omega \rangle = 0, \ \langle \omega, \bar{\omega}\rangle > 0\}
\end{equation}
which is a disjoint union of two copies of the 4-dimensional bounded symmetric domain of type IV.
We can consider ${\calD}(M)$ as a subdomain of ${\calD}(L_-)$ under the embedding $M \subset L_-$. 
We denote by $\Gamma_M$ the orthogonal group ${\rm O}(M)$ of $M$ and by $\tilde{\Gamma}_M$ the kernel of the map
$${\rm O}(M) \to {\rm O}(q_M).$$
Then $\Gamma_M/\tilde{\Gamma}_M \cong {\rm O}(q_M) \cong {\mathfrak S}_5 \times \{\pm 1\}$.
The quotient ${\calD}(M)/\tilde{\Gamma}_M$ is the moduli space of
{\it marked} Hessian quartic surfaces.

For a vector $r \in M^*$ with $r^2 < 0$,
we also define
$$r^{\perp} = \{ [\omega] \in {\calD}(M) \ | \ \langle \omega, r \rangle = 0\}.$$
Let $a \in A_M$ and let $m$ be a negative rational number. 
Define
$${\calH}_{a,m} = \sum_r \ r^{\perp}$$
where $r$ moves over the set of all vectors in $M^*$ satisfying $r\ {\rm mod}\ M = a$ and $r^2 = m$.  
We call ${\calH}_{a,m}$ the {\it Heegner divisor} of type $a$ and $m$.

\subsection{Proposition}\label{heegner1}
{\it 
A generic point of the Heegner divisor ${\calH}_{a, m}$ corresponds to the period of the Hessian quartic surfaces of the following cubic surfaces $S$.
If $a=0$ and $m=-2$, then $S$ has a node.  If $q_M(a)=1$ and $m=-1$, then $S$ has an Eckardt point.  
If $q_M(a)=m=-1/3$, then $S$ has no Sylvester forms.
}
\begin{proof}
The assertion follows from Dardanelli, van Geemen \cite{DvG}, Lemmas 2.2, 3.1, 5.1.
Also we have seen the assertion for $q_M(a)=1$ in Lemma \ref{norm1}.
\end{proof}

It follows from Sterk \cite{St}, Corollary 3.3 that any two vectors $r, s$ in $M$ satisfying
$$r^2 = s^2, \ \langle r, M\rangle = \langle s, M\rangle =: p{\bf Z}, \ (p > 0), \ r/p\ {\rm mod}\ M = s/p\ {\rm mod}\ M$$
are equivalent under the action of $\tilde{\Gamma}_M$.  In particular the image of each Heegner divisor
$${\calH}_{0, -2}, \  {\calH}_{a, -1} \ (a\in A_M, q_M(a) = 1), \ {\calH}_{a, -1/3} \ (a\in A_M, q_M(a) = -1/3)$$
in ${\calD}(M)/\tilde{\Gamma}_M$ is irreducible.
We denote by ${\calD}(M)^o$ the complement of the Heegner divisors ${\calH}_{0, -2}$ and
${\calH}_{a, -1/3}, \ a\in A_M, q_M(a) = -1/3$ in ${\calD}(M)$.
Let
\begin{equation}\label{moduli}
\Lambda = \{ \lambda \in {\bf P}^4 \ | \ \Delta_{\rm sing}(\lambda) \not= 0, \quad \lambda_i \not=0, \ i=1,...,5\}.
\end{equation} 
The symmetric group ${\mathfrak S}_5$ of degree 5 acts on $\Lambda$ as permutations of the coordinate of ${\bf P}^4$, 
and on ${\calD}(M)^o$ as the action of 
$\Gamma_M/\tilde{\Gamma}_M \cong {\mathfrak S}_5 \times \{\pm 1\}$.
Then the global Torelli type theorem for $K3$ surfaces and Proposition \ref{heegner1} imply the following proposition.
\subsection{Proposition}\label{torelli}{\rm (\cite{Koi}, Theorem 2.1)}
{\it The period map gives an ${\mathfrak S}_5$-equivariant embedding from $\Lambda$ to ${\calD}(M)^o/\tilde{\Gamma}_M$.}

\section{Automorphic forms}\label{forms}

In \cite{Kon}, the author constructed automorphic forms of weight 4 on ${\calD}(L_-)$ with respect to the group $\tilde{\rm O}(L_-)$
by using the theory of automorphic forms due to Borcherds \cite{B}.   
We recall this briefly.  

First we recall that $A_{L_-}$ is isomorphic to the orthogonal direct sum of 5 copies of $u$.
 A 5-dimensional subspace $V$ of $A_{L_-}$ is called a {\it maximal totally singular subspace} if
$V$ is generated by mutually orthogonal non-isotropic vectors $a_1,...,a_5$.  By using Borcherds theory \cite{B}, for each $V$, we can
associate a holomorphic automorphic form $F_V$ on ${\calD}(L_-)$ of weight 4 with respect to $\tilde{\rm O}(L_-)$.
Moreover the zero divisor of $F_V$ is 
\begin{equation}\label{period2}
\sum_{a \in V, q_{L_-}(a) = 1} \tilde{\calH}_a. 
\end{equation}
The linear system of these automorphic forms together with an another automorphic form of the same weight 
define an ${\rm O}^+(10, {\bf F})$-equivariant morphism from ${\calD}(L_-)/\tilde{\rm O}(L_-)$ to
${\bf P}^{186}$ which is birational onto its image
(In \cite{Kon}, there was a mistake pointed out and correced by Freitag and Manni \cite{FS}, Theorem 11.2).
There are $3^3\cdot 5\cdot 17\cdot 31$ maximal totally singular subspaces in $A_{L_-}$.

In the following, we consider the {\it restriction} of automorphic forms $F_V$ to
${\calD}(M)$.  Recall that $M$ is the orthogonal complement of $R \cong E_6(2)$ in $L_-$ (Lemma \ref{root-inv}).
Denote by $q_2$ the restriction of a discriminant quadratic form $q$ to the 2-Sylow subgroup.
Then 
$$(q_R)_2 \cong u\oplus u \oplus v, \quad (q_M)_2 \cong u \oplus v.$$
By using the relation $u \oplus u \cong v\oplus v$, we can see that
$$q_{L_-} \cong (q_R)_2 \oplus (q_M)_2.$$
Let $a_1, a_2, a_3$ be mutually orthogonal vectors with $q_R(a_i) =1$ $(i=1,2,3)$ in $A_R$.  
Then $a_1, a_2, a_3$ are mutually orthogonal non-isotropic vectors in $A_{L_-}$.
Let $b_1, b_2$ be mutually orthogonal vectors with $q_M(b_i)=1$ $(i=1,2)$ in $M^*/M$.
Then 
$$a_1, a_2, a_3, b_1, b_2$$ 
are mutually orthogonal non-isotropic vectors
in $A_{L_-}$, and hence they generate a maximal totally singular subspace $V$ in $A_{L_-}$.
There are 15 pairs $\{ b_1, b_2\}$ of mutually orthogonal non-isotropic vectors in $M^*/M$.
Thus, for fixed $a_1, a_2, a_3$, we have 15 maximal totally singular subspaces in $A_{L_-}$ of this type.

Now we study the restriction of Heegner divisors appeared in (\ref{period2}).
Let $V$ be a maximal totally singular subspace generated by
$a_1,a_2,a_3, b_1,b_2$.  Then $V$ contains 16 non-isotropic vectors
$$a_1, a_2, a_3, a_1+a_2+a_3, b_1, b_2, a_i+a_j+ b_k \ (1\leq i < j \leq 3, \ k = 1,2),$$
$$a_i + b_1 + b_2 \ (i=1,2,3), \ a_1+a_2+a_3+b_1+b_2.$$

\noindent
Obviously $\tilde{\calH}_a$ $(a =a_1, a_2,a_3, a_1+a_2+a_3)$
vanishes along ${\calD}(M)$.   On the other hand, if $r$ is a $(-4)$-vector in $L_-$ with $r/2 \ {\rm mod}\ L_- = b_j$, then $r \in M$ or the projection of $r$
into $M^*$ has a non-negative norm because the maximal norm of non-zero vectors in $R$ is $-4$.  
In the later case, the hyperplane $r^{\perp}$ does not meet with ${\calD}(M)$.  Therefore 
$\tilde{\calH}_{b_j}$ $(j=1,2)$ cuts the Heegner divisor
${\calH}_{b_j, -1}$ on ${\calD}(M)$. 
In case $a= a_i+a_j+b_k$, $a_i+a_j$ is represented by $r/2$ with $r\in R$,
$r^2 = -8m$, $m$ a positive integer, because $a_i+a_j$ is non-zero isotropic vector and 
the maximal norm of non-zero vectors in $R$ is $-4$.
This implies that $b_k$ is represented by a positive norm vector in $M$, 
and hence $\tilde{\calH}_a$ does not intersect with ${\calD}(M)$ and its boundary.
Similary in case $a= a_i + b_1+b_2$ or $a=a_1+a_2+a_3+b_1+b_2$, 
$a_i$ and $a_1+a_2+a_3$ are represented by $r/2$ with $r \in R$, $r^2 = -4m$ where $m$ is a positive integer.  This implies that
$b_1+b_2$ is represented by a non-zero isotropic vector or a positive norm vector in $M$, and hence
$\tilde{\calH}_a$ does not meet with the interior of ${\calD}(M)$.

Now recall that $E_6$ contains 72 roots (\cite{Bou}, Planche V), and hence $R=E_6(2)$ contains 72 $(-4)$-vectors.  On the other hand, the number of non-isotropic
vectors of the quadratic form $(q_R)_2 = u_2\oplus u_2\oplus v_2$ is 36.   
By sending each $(-4)$-vector $\pm r$ in $R$ to $r/2 \ {\rm mod}\ R$ in $(q_R)_2$, 
we have a bijective correspondence between the set of $(-4)$-vectors in $R$ modulo $\pm 1$ and the set of non-isotropic vectors in $(q_R)_2$.
Thus each $F_V$ vanishes along 4 hyperplanes $(\pm r)^{\perp}$ where $r \in R$ with $r^2 =-4$ 
corresponding to 4 non-isotropic vectors $a_1, a_2, a_3, a_1+a_2+a_3$.
By a method given in \cite{B1}, that is,  by first dividing $F_V$ by a product of linear forms vanishing on the divisors associated to the four
$(-4)$-vectors and restricting it to ${\calD}(M)$ we can get an automorphic form $F_V|M$ on ${\calD}(M)$ with respect to the group
$\tilde{\Gamma}_M$.  Then the weight of $F_V|M$ is the weight of
$F_V$ plus half the number of $(-4)$-vectors in $R$ corresponding to
$a_1, a_2, a_3, a_1+a_2+a_3$, that is, $4+4=8$.
We now conclude

\subsection{Theorem}\label{rest}
{\it Let $V$ be a maximal totally singular subspace generated by
$\{a_1,a_2,a_3, b_1,b_2\}$.  Then $F_V|M$ is a holomorphic automorphic form on ${\calD}(M)$ 
of weight {\rm 8} with respect to the group $\tilde{\Gamma}_M$ whose zero divisor is
${\calH}_{b_1,-1} + {\calH}_{b_2,-1}.$ }

\medskip

As mentioned in Proposition \ref{heegner1}, a generic point of the Heegner divisor ${\calH}_{b_i, -1}$ is the period of
the Hessian quartic surface of a cubic surface with an Eckardt point.  Assume that non-isotropic vector $b_1$ or $b_2$ corresponds to the condition
$\lambda_i = \lambda_j$ or $\lambda_k= \lambda_l$ respectively.  Then the condition that $b_1$ is orthogonal to $b_2$ is equivalent
to that all $i,j,k,l$ are different (see Lemma \ref{norm1}).  By identifying $\Lambda$ and an open subset in ${\calD}(M)$ (see Proposition \ref{torelli}),
we have the following theorem.

\subsection{Theorem}\label{rest2}
{\it As divisors on $\Lambda$,}
$$(F_V|M) =  ((\lambda_i -\lambda_j)(\lambda_k -\lambda_l)).$$

\subsection{Remark}\label{}
{\it We can easily see that fifteen $(\lambda_i - \lambda_j)(\lambda_k - \lambda_l)$, where $i, j, k, l$ are distinct, generate a $5$-dimensional space $W$ 
of quadrics on ${\bf P}^4$ whose base locus is the sum of
five lines defined by $\lambda_i = \lambda_j = \lambda_k$ $(1\leq i < j < k\leq 5)$ meeting at $(1:1:1:1:1)$.
Let
$$\xi = (\lambda_1 -\lambda_2)(\lambda_3-\lambda_5), \ 
\eta = (\lambda_1 -\lambda_3)(\lambda_4-\lambda_5),\ 
\zeta = (\lambda_1 -\lambda_4)(\lambda_2-\lambda_5),$$
$$\xi' =  (\lambda_1 -\lambda_2)(\lambda_4-\lambda_5),\ 
\eta'= (\lambda_1 -\lambda_3)(\lambda_2-\lambda_5), \ 
 \zeta'= (\lambda_1 -\lambda_4)(\lambda_3-\lambda_5).$$
The $\xi, \xi', \eta, \eta', \zeta, \zeta'$ generate $W$ and satisfy the relations
$$\xi + \eta + \zeta = \xi' + \eta' + \zeta', \ \xi \eta \zeta = \xi' \eta' \zeta'.$$
These relations define the Segre cubic $3$-fold $\calS_3$ {\rm (}Baker {\rm \cite{Ba}}, Hunt {\rm \cite{H}, \S 3.2.2)}.  
Since the restriction of $W$ on a general hyperplane ${\bf P}^3$ in ${\bf P}^4$ is the linear system of
quadrics with five points as its base locus, it gives a birational map from ${\bf P}^3$ to $\calS_3$ {\rm (}Hunt  {\rm \cite{H}, Theorem 3.2.1)}.  
Thus the linear system $W$ gives a dominant rational map
from ${\bf P}^4$ to $\calS_3$.  
}

\subsection{Remark}\label{}
{\it Borcherds {\rm \cite{B0}} constructed an automorphic form $\Phi_4$ on ${\calD}(L_-)$ of weight $4$ whose zero divisor is the Heegner divisor 
$\tilde{\calH}$
associated to $(-2)$-vectors in $L_-$.  Since $R=E_6(2)$ has no $(-2)$-vectors, the restriction of $\Phi_4$ defines an automorphic form on ${\calD}(M)$ of weight $4$.  Let $r$ be a $(-2)$-vector in $L_-$ and let
$$r = r_1 + r_2, \ r_1 \in R^*, \ r_2 \in M^*.$$
Assume that $r_1\not= 0$.  Since $R$ contains no $(-2)$-vectors, $r_2 \not=0$.
Since $M \oplus R$ has index $3$ in $L_-$, 
$3r_1 \in R, \ 3r_2 \in M$.  By Lemma {\rm \ref{discri}}, $q_M(r_2\ {\rm mod}\ M) = -4/3$.  Hence $q_{R}(r_1 \ {\rm mod}\ R) = -2/3$.
Since $R$ contains no $(-6)$-vectors, $(r_1)^2 \leq -8/3$ and hence $(r_2)^2 > 0$.
Therefore $r^{\perp}$ does not intersect with ${\calD}(M)$.  This implies that 
if the projection of a $(-2)$-vector in $L_-$ into $M^*$  has a negative norm, then it is a $(-2)$-vector.
Thus the restriction
of $\Phi_4$ is an automorphic form on ${\calD}(M)$ of weight $4$ whose zero divisor is the Heegner divisor ${\calH}_{0,-2}$ 
associated to $(-2)$-vectors in $M$.
The corresponding cubic surfaces are nodal {\rm (}Proposition {\rm \ref{heegner1})}.
}

\end{document}